# Critical scaling of stochastic epidemic models[*]

**Steven P. Lalley**[1]

*University of Chicago*

**To Piet Groeneboom, on the occasion of his 39th birthday.**

**Abstract:** In the simple mean-field *SIS* and *SIR* epidemic models, infection is transmitted from infectious to susceptible members of a finite population by independent $p-$coin tosses. Spatial variants of these models are proposed, in which finite populations of size $N$ are situated at the sites of a lattice and infectious contacts are limited to individuals at neighboring sites. Scaling laws for both the mean-field and spatial models are given when the infection parameter $p$ is such that the epidemics are *critical*. It is shown that in all cases there is a *critical threshold* for the numbers initially infected: below the threshold, the epidemic evolves in essentially the same manner as its *branching envelope*, but at the threshold evolves like a branching process with a size-dependent drift.

## 1. Stochastic epidemic models

### 1.1. Mean-field models

The simplest and most thoroughly studied stochastic models of epidemics are *mean-field* models, in which all individuals of a finite population interact in the same manner. In these models, a contagious disease is transmitted among individuals of a homogeneous population of size $N$. In the simple *SIS epidemic*, individuals are at any time either *infected* or *susceptible*; infected individuals remain infected for one unit of time and then become susceptible. In the simple *SIR epidemic* (more commonly known as the *Reed-Frost* model), individuals are either *infected, susceptible,* or *recovered*; infected individuals remain infected for one unit of time, after which they recover and acquire permanent immunity from future infection. In both models, the mechanism by which infection occurs is random: At each time, for any pair $(i,s)$ of an infected and a susceptible individual, the disease is transmitted from $i$ to $s$ with probability $p = p_N$. These transmission events are mutually independent. Thus, in both the *SIR* and the *SIS* model, the number $J_{t+1} = J_t^N$ of infected individuals at time $t+1$ is given by

$$J_{t+1} = \sum_{s=1}^{S_t} \xi_s, \tag{1}$$

[*]Supported by NSF Grant DMS-04-05102.
[1]University of Chicago, Department of Statistics, 5734 S. University Avenue, Eckhart 118, Chicago, Illinois 60637, USA, e-mail: lalley@galton.uchicago.edu

where $S_t = S_t^N$ is the number of susceptibles at time $t$ and the random variables $\xi_s$ are, conditional on the history of the epidemic to time $t$, independent, identically distributed Bernoulli-$1 - (1-p)^{J_t}$. In the $SIR$ model,

$$R_{t+1} = R_t + J_t \quad \text{and} \tag{2}$$
$$S_{t+1} = S_t - J_{t+1},$$

where $R_t$ is the number of recovered individuals at time $t$, while in the $SIS$ model,

$$S_{t+1} = S_t + J_t - J_{t+1}. \tag{3}$$

In either model, the epidemic ends at the first time $T$ when $J_T = 0$. The most basic and interesting questions concerning these models have to do with the *duration* $T$ and *size* $\sum_{t \leq T} J_t$ of the epidemic and their dependence on the infection parameter $p_N$ and the initial conditions.

### 1.2. Spatial SIR and SIS epidemics

In the simple $SIS$ and $SIR$ epidemics, no allowance is made for geographic or social stratifications of the population, nor for variability in susceptibility or degree of contagiousness. Following are descriptions of simple stochastic models that incorporate a geographic stratification of a population. We shall call these the (spatial) $SIS-d$ and $SIR-d$ epidemics, with $d$ denoting the spatial dimension.

Assume that at each lattice point $x \in \mathbb{Z}^d$ is a homogeneous population of $N_x$ individuals, each of whom may at any time be either susceptible or infected, or (in the $SIR$ variants) recovered. These populations may be thought of as "villages". As in the mean-field models, infected individuals remain contagious for one unit of time, after which they recover with immunity from future infection (in the $SIR$ variants) or once again become susceptible (in the $SIS$ models). At each time $t = 0, 1, 2, \ldots$, for each pair $(i_x, s_y)$ of an infected individual located at $x$ and a susceptible individual at $y$, the disease spreads from $i_x$ to $s_y$ with probability $\alpha(x, y)$.

The simple Reed-Frost and stochastic logistic epidemics described in section 1.1 terminate with probability one, regardless of the value of the infection parameter $p$, because the population is finite. For the spatial $SIS$ and $SIR$ models this is no longer necessarily the case: If $\sum_{x \in \mathbb{Z}^d} N_x = \infty$ then, depending on the value of the parameter $p$ and the dimension $d$, the epidemic may persist forever with positive probability. (For instance, if $N_x = 1$ for all $x$ and $\alpha(x, y) = p$ for nearest neighbor pairs $x, y$ but $\alpha(x, y) = 0$ otherwise, then the $SIS-d$ epidemic is just oriented percolation on $\mathbb{Z}^{d+1}$, which is known to survive with positive probability if $p$ exceeds a critical value $p_c < 1$ [6].) Obvious questions of interest center on how the epidemic spreads through space, and in cases where it eventually dies out, how far it spreads.

The figure below shows a simulation of an $SIS$-1 epidemic with village size $N_x = 20224$ and infection parameter $1/60672$. At time 0 there were 2048 infected individuals at site 0; all other individuals were healthy. The epidemic lasted 713 generations (only the first 450 are shown).

### 1.3. Epidemic models and random graphs

All of the models described above have equivalent descriptions as structured random graphs, that is, percolation processes. Consider for definiteness the simple $SIR$ (Reed-Frost) epidemic. In this model, no individual may be infected more than once; furthermore, for any pair $x, y$ of individuals, there will be at most one opportunity for infection to pass from $x$ to $y$ or from $y$ to $x$ during the course of the epidemic. Thus, one could simulate the epidemic by first tossing a $p-$coin for every pair $x, y$, drawing an edge between $x$ and $y$ for each coin toss resulting in a Head, and then using the resulting (Erdös-Renyi) random graph determined by these edges to determine the course of infection in the epidemic. In detail: If $Y_0$ is the set of infected individuals at time 0, then the set $Y_1$ of individuals infected at time 1 consists of all $x \notin Y_0$ that are connected to individuals in $Y_0$, and for any subsequent time $n$, the set $Y_{n+1}$ of individuals infected at time $n+1$ consists of all $x \notin \cup_{j \leq n} Y_j$ who are connected to individuals in $Y_n$. Note that the set of individuals ultimately infected during the course of the epidemic is the union of those connected components of the random graph containing at least one vertex in $Y_0$.

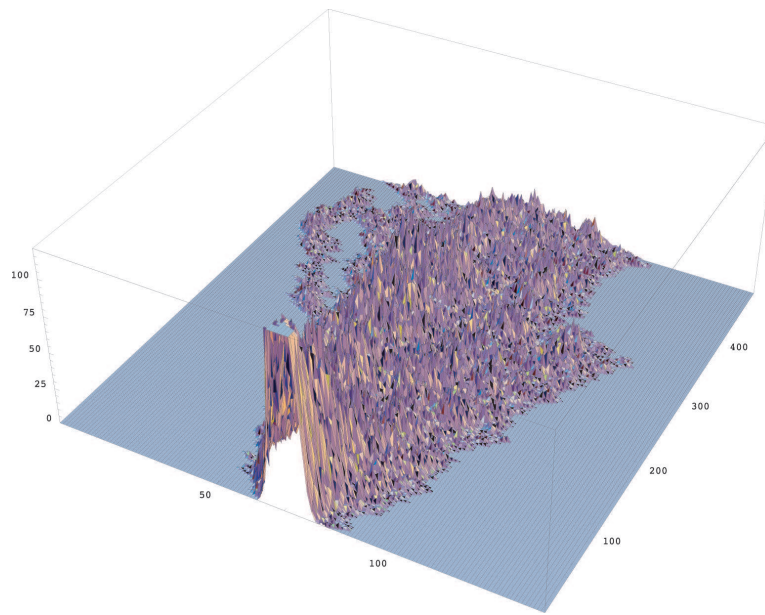

Fig 1.

Similar random graph descriptions may be given for the simple *SIS* and the spatial *SIS* and *SIR* epidemic models.

### *1.4. Branching envelopes of epidemics*

For each of the stochastic epidemic models discussed above there is an associated branching process that serves, in a certain sense, as a "tangent" to the epidemic. We shall refer to this branching process as the *branching envelope* of the epidemic. The branching envelopes of the simple mean-field epidemics are ordinary Galton-Watson processes; the envelopes of the spatial epidemics are branching random walks. There is a natural coupling of each epidemic with its branching envelope in which the set of infected individuals in the epidemic is at each time (and in the spatial models, at each location) dominated by the corresponding set of individuals in the branching envelope.

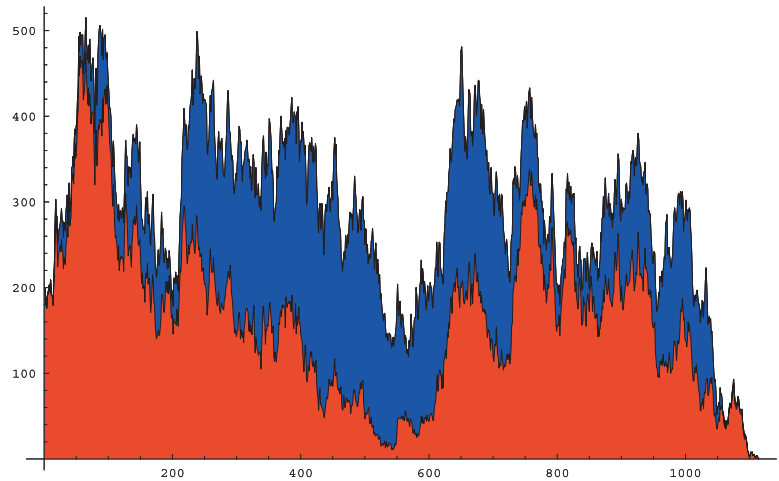

Fig 2.

Following is a detailed description of the natural coupling of the simple *SIS* epidemic with its branching envelope. The branching envelope is a Galton-Watson process $Z_n$ with offspring distribution Binomial-$(N, p)$, where $p$ is the infection parameter of the epidemic, and whose initial generation $Z_0$ coincides with the set of individuals who are infected at time 0. Particles in the Galton-Watson process are marked *red* or *blue*: red particles represent infected individuals in the coupled epidemic, while blue offspring of red parents represent attempted infections that are not allowed because the attempt is made on an individual who is not susceptible, or has already been infected by another contagious individual. Colors are assigned as follows: (1) Offspring of blue particles are always blue. (2) Each red particle reproduces by tossing a $p-$coin $N$ times, once for each individual $i$ in the population. Each Head counts as an offspring, and each represents an *attempted* infection. If several red particles attempt to infect the same individual $i$, exactly one of these is marked as a success (red), and the others are marked as failures (blue). Also, if an attempt is made to infect an individual who is not susceptible, the corresponding particle is colored blue. Clearly, the collection of all particles (red and blue) evolves as a Galton-Watson process, while the collection of red particles evolves as the infected set in the SIS epidemic. See the figure below for a typical evolution of the coupling in a population of size $N = 80,000$ with $p = 1/80000$ and 200 individuals initially infected.

## 2. Critical behavior: mean-field case

When studying the behavior of the simple *SIR* and *SIS* epidemics in large populations, it is natural to consider the scaling $p = p_N = \lambda_N/N$ for the infection parameter $p$. In this scaling, $\lambda = \lambda_N$ is the mean of the offspring distribution in the branching envelope. If $\lambda < 1$ then the epidemic will end quickly, even if a large number of individuals are initially infected. On the other hand, if $\lambda > 1$ then with positive probability (approximately one minus the extinction probability for the associated branching envelope), even if only one individual is initially infected, the epidemic will be large, with a positive limiting fraction of the population eventually being infected. The large-$N$ behavior of the size of the *SIR* epidemic in this case is well understood: see for example [12] and [14].

### 2.1. Critical scaling: size of the epidemic

The behavior of both the *SIS* and *SIR* epidemics is more interesting in the *critical case* $\lambda_N \approx 1$. When the set of individuals initially infected is sufficiently small relative to the population size, the epidemic can be expected to evolve in much the same manner as a critical Galton-Watson process with Poisson-1 offspring distribution. However, when the size of the initially infected set passes a certain *critical threshold*, then the epidemic will begin to deviate substantially from the branching envelope. For the *SIR* case, the critical threshold was (implicitly) shown by [11] and [1] (see also [12]) to be at $N^{1/3}$, and that the *critical scaling window* is of width $N^{-4/3}$:

**Theorem 1** ([11], [1]). *Assume that $p_N = 1/N + a/N^{4/3} + o(n^{-4/3})$, and that the number $J_0^N$ of initially infected individuals is such that $J_0^N/N^{1/3} \to b$ as the population size $N \to \infty$. Then as $N \to \infty$, the size $U_N := \sum_t J_t$ obeys the asymptotic law*

$$U_N/N^{2/3} \xrightarrow{\mathcal{D}} T_b \tag{4}$$

*where $T_b$ is the first passage time to the level $b$ by $W_t + t^2/2 + at$, and $W_t$ is a standard Wiener process.*

The distribution of the first passage time $T_b$ can be given in closed form: See [11], also [8], [13].

For the critical *SIS* epidemic, the critical threshold is at $N^{1/2}$, and the critical scaling window is of width $N^{-3/2}$:

**Theorem 2** ([4]). *Assume that $p_N = 1/N + a/N^{3/2} + o(n^{-3/2})$, and that the initial number of infected individuals satisfies $J_0^N \sim bN^{1/2}$ as $N \to \infty$. Then the total number of infections $U_N := \sum_t J_t$ during the course of the epidemic obeys*

$$U_N/N \xrightarrow{\mathcal{D}} \tau(b-a; -a) \tag{5}$$

*where $\tau(x; y)$ is the time of first passage to $y$ by a standard Ornstein-Uhlenbeck process started at $x$.*

### 2.2. Critical scaling: time evolution of the epidemic

For both the *SIR* and *SIS* epidemics, if the number of individuals initially infected is much below the critical threshold then the evolution of the epidemic will not differ noticeably from that of its branching envelope. It was observed by [7] (and proved by [9]) that a (near-) critical Galton-Watson process initiated by a large number $M$ of individuals behaves, after appropriate rescaling, approximately as a *Feller diffusion*: In particular, if $Z_n^M$ is the size of the $n$th generation of a Galton-Watson with $Z_0^M \sim bM$ with offspring distribution Poisson$(1 + a/M)$ then as $M \to \infty$,

$$Z_{[Mt]}^M/M \xrightarrow{\mathcal{D}} Y_t \tag{6}$$

where $Y_t$ satisfies the stochastic differential equation

$$dY_t = aY_t \, dt + \sqrt{Y_t} \, dW_t, \tag{7}$$
$$Y_0 = b.$$

What happens at the critical threshold, in both the *SIR* and *SIS* epidemics, is that the deviation from the branching envelope exhibits itself as a size-dependent drift in the limiting diffusion:

**Theorem 3** ([4]). *Let $J^N(n) = J^N_{[n]}$ be the number infected in the nth generation of a simple* SIS *epidemic in a population of size $N$. Then under the hypotheses of Theorem 2,*

$$J^N(\sqrt{N}t)/\sqrt{N} \xrightarrow{\mathcal{D}} Y_t \tag{8}$$

*where $Y_0 = b$ and $Y_t$ obeys the stochastic differential equation*

$$dY_t = (aY_t - Y_t^2)\,dt + \sqrt{Y_t}\,dW_t \tag{9}$$

Note that the diffusion (9) has an entrance boundary at $\infty$, so that it is possible to define a version $Y_t$ of the process with initial condition $Y_0 = 0$. When the *SIS* epidemic is begun with $J_0^N \gg \sqrt{N}$ initially infected, the number $J_t^N$ infected will rapidly drop (over the first $\varepsilon\sqrt{N}$ generations) until reaching a level of order $\sqrt{N}$, and then evolve as predicted by (8). The following figure depicts a typical evolution in a population of size $N = 80,000$, with infection parameter $p = 1/N$ and $I_0 = 10,000$ initially infected.

**Theorem 4** ([4]). *Let $J^N(n) = J^N_{[n]}$ and $R^N(n) = R^N_{[n]}$ be the numbers of infected and recovered individuals in the nth generation of a simple* SIR *epidemic in a population of size $N$. Then under the hypotheses of Theorem 1,*

$$\begin{pmatrix} N^{-1/3}J^N(N^{1/3}t) \\ N^{-2/3}R^N(N^{1/3}t) \end{pmatrix} \xrightarrow{\mathcal{D}} \begin{pmatrix} J(t) \\ R(t) \end{pmatrix} \tag{10}$$

*where $J_0 = b$, $R_0 = 0$, and*

$$\begin{aligned} dJ(t) &= (aJ(t) - J(t)R(t))\,dt + \sqrt{J(t)}\,dW_t, \\ dR(t) &= J(t)\,dt. \end{aligned} \tag{11}$$

Theorems 1–2 can be deduced from Theorems 3–4 by simple time-change arguments (see [4]).

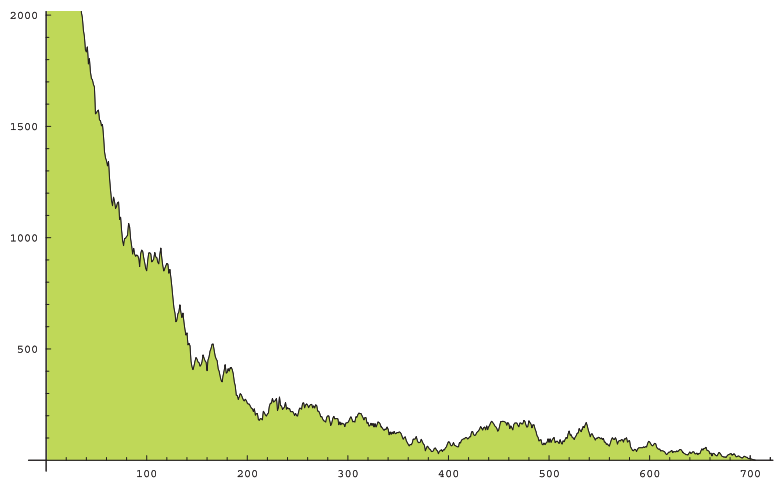

Fig 3.

## 2.3. Critical scaling: heuristics

The critical thresholds for the $SIS-d$ and $SIR-d$ epidemics can be guessed by simple comparison arguments using the standard couplings of the epidemics with their branching envelopes. Consider first the critical $SIS$ epidemic in a population of size $N$. Recall (Section 1.4) that the branching envelope is a critical Galton-Watson process whose offspring distribution is Binomial-$(N, 1/N)$. The particles of this Galton-Watson process are marked red or blue, in such a way that in each generation the number of red particles coincides with the number of infected individuals in the $SIS$ epidemic. Offspring of blue particles are always blue, but offspring of red particles may be either red or blue; the blue offspring of red parents in each generation represent attempted infections that are suppressed.

Assume that initially there are $N^\alpha$ infected individuals and thus also $N^\alpha$ individuals in the zeroth generation of the branching envelope. By Feller's theorem, we may expect that the extinction time of the branching envelope will be on the order $N^\alpha$, and that in each generation up to (shortly before) extinction the branching process will have order $N^\alpha$ individuals. If $\alpha$ is small enough that the SIS epidemic obeys the same rough asymptotics (that is, stays alive for $O(N^\alpha)$ generations and has $O(N^\alpha)$ infected individuals in each generation), then the number of blue offspring of red parents in each generation will be on the order $N \times (N^{2\alpha}/N^2)$ (because for each of the $N$ individuals of the population, the chance of a double infection is about $N^{2\alpha}/N^2$). Since the duration of the epidemic will be of the same rough order of magnitude as the size of the infected set in each generation, there should be at most $O(1)$ blue offspring of red parents in any generation (if there were more, the red population would die out long before the blue). Thus, the critical threshold must be at $N^{1/2}$.

A similar argument applies for the $SIR$ epidemic. The branching envelope of the critical $SIR$ is once again a critical Galton-Watson process with offspring distribution Binomial-$(N, 1/N)$, with constituent particles again labeled red or blue, red particles representing infected individuals in the epidemic. The rule by which red particles reproduce is as follows: Each red particle tosses a $p-$coin $N$ times once for each individual $i$ in the population. Each Head counts as an offspring, and represents an attempted infection. However, if a Head occurs on a toss at individual $i$ where $i$ was infected in an earlier generation, then the Head results in a *blue* offspring. Similarly, if more than one red particle tosses a Head at an individual $i$ which has not been infected earlier, then one of these is labeled red and the excess are all labeled blue.

Assume that initially there are $N^\alpha$ infected individuals. As before, we may expect that the extinction time of the branching envelope will be on the order $N^\alpha$, and that in each generation up to extinction the branching process will have order $N^\alpha$ individuals. If $\alpha$ is small enough, the extinction time and the size of the red population will also be $O(N^\alpha)$. Consequently, the size of the *recovered* population will be (for all but the first few generations) on order $N^{2\alpha}$. Thus, in each generation, the number of blue offspring of red parents will be on order $(N^{2\alpha}/N) \times N^\alpha$ (because the chance that a recovered individual is chosen for attempted infection by an infected individual is $O(N^\alpha/N)$). Therefore, by similar reasoning as in the $SIS$ case, the critical threshold is at $N^{1/3}$, as this is where the the number of blue offspring of red parents in each generation is $O(1)$.

## 3. Critical behavior: *SIS*-1 and *SIR*-1 Spatial epidemics

Consider now the spatial *SIS-d* and *SIR-d* epidemic models on the $d$-dimensional integer lattice $\mathbb{Z}^d$. Assume that the village size $N_x = N$ is the same for all sites $x \in \mathbb{Z}^d$, and that the infection probabilities $\alpha(x, y)$ are nearest neighbor homogeneous, and uniform, that is,

$$\alpha(x, y) = \begin{cases} p = p_N, & \text{if } |x - y| \leq 1; \\ 0, & \text{otherwise.} \end{cases} \tag{12}$$

### 3.1. Scaling limits of branching random walks

The *branching envelope* of a spatial $SIS-d$ or $SIR-d$ epidemic is a nearest neighbor *branching random walk* on the integer lattice $\mathbb{Z}^d$. This evolves as follows: Any particle located at site $x$ at time $t$ lives for one unit of time and then reproduces, placing random numbers $\xi_y$ of offspring at the sites $y$ such that $|y - x| \leq 1$. The random variables $\xi_y$ are mutually independent, with Binomial-$(N, p_N)$ distributions.

The analogue for branching random walks of Feller's theorem for Galton-Watson processes is *Watanabe's theorem*. This asserts that, after suitable rescaling, as the particle density increases, critical branching random walks converge to a limit, the *Dawson-Watanabe* process, also known as *super Brownian motion*. A precise statement follows: Consider a sequence of branching random walks, indexed by $M = 1, 2, \ldots$, with offspring distribution Binomial-$(N, p_M)$ as above, and

$$p_M = p_{N,M} = \frac{1}{(2d + 1)N} - \frac{a}{NM}. \tag{13}$$

(Note: $N$ may depend on $M$.) The rescaled *measure-valued process* $X_t^M$ associated with the $M$th branching random walk puts mass $1/M$ at location $x/\sqrt{M}$ and time $t$ for each particle of the branching random walk that is located at site $x$ at time $[Mt]$. (Note: The branching random walk is a discrete-time process, but the associated measure-valued process runs in continuous time.)

**Watanabe's theorem** ([15]). *Assume that the initial values $X_0^M$ converge weakly (as finite Borel measures on $\mathbb{R}^d$) to a limit measure $X_0$. Then under the hypothesis (13) the measure-valued processes $X_t^M$ converge in law as $M \to \infty$ to a limit $X_t$:*

$$X_t^M \Longrightarrow X_t. \tag{14}$$

The limit process is the *Dawson-Watanabe* process with killing rate $a$ and initial value $X_0$. (The term *killing rate* is used because the process can be obtained from the "standard" Dawson-Watanabe process ($a = 0$) by elimination of mass at constant rate $a$.) The Dawson-Watanabe process $X_t$ with killing rate $a$ can be characterized by a martingale property: For each test function $\phi \in C_c^2(\mathbb{R}^d)$,

$$\langle X_t, \phi \rangle - \langle X_0, \phi \rangle - \frac{\sigma}{2} \int_0^t \langle X_s, \Delta \phi \rangle \, ds + a \int_0^t \langle X_s, \varphi \rangle \, ds \tag{15}$$

is a martingale. Here $\sigma^2 = 2d/(2d + 1)$ is the variance of the random walk kernel naturally associated with the branching random walks. It is known [10] that in dimension $d = 1$ the random measure $X_t$ is for each $t$ absolutely continuous relative to Lebesgue measure, and the Radon-Nikodym derivative $X(t, x)$ is jointly continuous in $t, x$ (for $t > 0$). In dimensions $d \geq 2$ the measure $X_t$ is almost surely singular, and is supported by a Borel set of Hausdorff dimension 2 [3].

### 3.2. Spatial SIS-1 and SIR-1 epidemics: critical scaling

As in the mean-field case, there are critical thresholds for the *SIS*-1 and *SIR*-1 epidemics at which they begin to deviate noticeably from their branching envelopes. These are at $N^{2/3}$ and $N^{2/5}$, respectively:

**Theorem 5** ([5]). *Fix $\alpha > 0$, and let $X_t^N$ be the discrete-time measure-valued process obtained from an* SIS$-1$ *or an* SIR$-1$ *epidemic on a one-dimensional grid of size-N villages by attaching mass $1/N^\alpha$ to the point $(t, x/N^{\alpha/2})$ for each infected individual at site $x$ at time $[tN^\alpha]$. Assume that $X_0^N$ converges weakly to a limit measure $X_0$ as the village size $N \to \infty$. Then as $N \to \infty$,*

$$(16) \qquad X_{[N^\alpha t]}^N \xrightarrow{\mathcal{D}} X_t,$$

*where $X_t$ is a measure-valued process with initial value $X_0$ whose law depends on the value of $\alpha$ and the type of epidemic (SIS or SIR) as follows:*

(a) SIS: If $\alpha < \frac{2}{3}$ then $X_t$ is a Dawson-Watanabe process with variance $\sigma^2$.
(b) SIS: If $\alpha = \frac{2}{3}$ then $X_t$ is a Dawson-Watanabe process with variance $\sigma^2$ and killing rate

$$(17) \qquad \theta(x,t) = X(x,t)/2.$$

(c) SIR: If $\alpha < \frac{2}{5}$ then $X_t$ is a Dawson-Watanabe process with variance $\sigma^2$.
(d) SIR: If $\alpha = \frac{2}{5}$ then $X_t$ is a Dawson-Watanabe process with variance $\sigma^2$ and killing rate

$$(18) \qquad \theta(x,t) = X(x,t) \int_0^t X(x,s)\,ds.$$

The Dawson-Watanabe process with variance $\sigma^2$ and (continuous, adapted) killing rate $\theta(t, x, \omega)$ is characterized [2] by a martingale problem similar to (15) above: for each test function $\phi \in C_c^2(\mathbb{R})$,

$$(19) \qquad \langle X_t, \phi \rangle - \langle X_0, \phi \rangle - \frac{\sigma}{2}\int_0^t \langle X_s, \Delta\phi \rangle\,ds + \int_0^t \langle X_s, \theta\varphi \rangle\,ds$$

is a martingale. The law of this process is mutually absolutely continuous relative to that of the Dawson-Watanabe process with no killing, and there is an explicit formula for the Radon-Nikodym derivative – see [2].

### 3.3. Critical scaling for spatial epidemics: heuristics

Arguments similar to those given above for the mean-field *SIS* and *SIR* epidemics can be used to guess the critical thresholds for the spatial *SIS*-1 and *SIR*-1 epidemics. For the spatial epidemics, the associated branching envelopes are branching random walks. In the standard couplings, particles of the branching envelope are labeled either red or blue, with the red particles representing infected individuals in the epidemics. As in the mean-field cases, offspring of blue particles are always blue, but offspring of red particles may be either red or blue; blue offspring of red parents represent attempted infections that are suppressed. These may be viewed as an *attrition* of the red population (since blue particles created by red parents are potential red offspring that are not realized!).

Consider first the *SIS*-1 epidemic. Assume that initially there are $N^\alpha$ particles, distributed (say) uniformly among the $N^{\alpha/2}$ sites nearest the origin. Then by Feller's limit theorem (recall that the total population size in a branching random walk is a Galton-Watson process), the branching envelope can be expected to survive for $O_P(N^\alpha)$ generations, and at any time prior to extinction the population will have $O_P(N^\alpha)$ members. These will be distributed among the sites at distance $O_P(N^{\alpha/2})$ from the origin, and therefore in dimension $d = 1$ there should be about $O_P(N^{\alpha/2})$ particles per site. Consequently, for the $SIS-1$ epidemic, the rate of attrition per site per generation should be $O_P(N^{\alpha-1})$, and so the total attrition rate per generation should be $O_P(N^{3\alpha/2-1})$. If $\alpha = 2/3$, then the total attrition rate per generation will be $O_P(1)$, just enough so that the total attrition through the duration of the branching random walk envelope will be on the same order of magnitude as the population size $N^\alpha$.

For the $SIR-1$ epidemic there is a similar heuristic calculation. As for the $SIS-1$ epidemic, the branching envelope will survive for $O_P(N^\alpha)$ generations, and up to the time of extinction the population should have $O_P(N^\alpha)$ individuals, about $O_P(N^{\alpha/2})$ per site. Therefore, through $N^\alpha$ generations, about $N^\alpha \times N^{\alpha/2}$ numbers $j$ will be retired, and so the attrition rate per site per generation should be $O_P(N^{\alpha/2} \times N^{3\alpha/2})$, making the total attrition rate per generation $O_P(N^{5\alpha/2})$. Hence, if $\alpha = 2/5$ then the total attrition per generation should be $O_P(1)$, just enough so that the total attrition through the duration of the branching random walk envelope will be on the same order of magnitude as the population size.

### 3.4. Critical scaling in dimensions $d \geq 2$

In higher dimensions, the critical behavior of the *SIS-d* and *SIR-d* epidemics appears to be considerably different. We expect that there will be no analogous threshold effect, in particular, we expect that the epidemic will behave in the same manner as the branching envelope up to the point where the infected set is a positive fraction of the total population. This is because in dimensions $d \geq 2$, the particles of a critical branching random walk quickly spread out, so that (after a short initial period) there are only $O_P(1)$ (in dimensions $d \geq 3$) or $O_P(\log N)$ (in dimension $d = 2$) particles per site. (With $N$ particles initially, a critical branching random walk typically lives $O(N)$ generations, and particles are distributed among the sites at distance $O(\sqrt{N})$ from the origin; in dimensions $d \geq 2$, there are $O(N^{d/2})$ such sites, enough to accomodate the $O(N)$ particles of the branching random walk without crowding.) Consequently, the rate at which "multiple" infections are attempted (that is, attempts by more than one contagious individual to simultaneously infect the same susceptible) is only of order $O_P(1/N)$ (or, in dimension $d = 2$, order $O_P(\log N/N)$).

The interesting questions regarding the evolution of critical epidemics in dimensions $d \geq 2$ center on the initial stages, in the relatively small amount of time (order $o(N)$ generations) in which the particles spread out from their initial sites. These will be discussed in the forthcoming University of Chicago Ph. D. dissertation of Xinghua Zheng.

### 3.5. Weak convergence of densities

There is an obvious gap in the heuristic argument of Section 3.3 above: Even if the total number of infected individuals is, as expected, on the order $N^\alpha$, and even if these are concentrated in the sites at distance on the order $N^{\alpha/2}$ from the origin, it is by no means obvious that these will distribute themselves uniformly (or at least *locally* uniformly) among these sites. The key step in filling this gap in the argument is to show that the particles of the branching envelope distribute themselves more or less uniformly on scales smaller than $N^{\alpha/2}$.

Consider, as in Section 3.1, a sequence of branching random walks, indexed by $M = 1, 2, \ldots$, with offspring distribution Binomial-$(N, p_M)$ as above, and $p_M$ given by (13). Let $Y_t^M(x)$ be the number of particles at site $x$ at time $[t]$, and let $X^M(t, x)$ be the continuous function of $t \geq 0$ and $x \in \mathbb{R}$ obtained by linear interpolation from the values

$$(20) \qquad X^M(t, x) = \frac{Y_{Mt}(\sqrt{M}x)}{\sqrt{M}} \quad \text{for} \quad Mt \in \mathbb{Z}_+ \text{ and } \sqrt{M}x \in \mathbb{Z}.$$

**Theorem 6** ([5]). *Assume that $d = 1$. Assume also that the initial particle configuration is such that all particles are located in an interval $[-\kappa\sqrt{M}, \kappa\sqrt{M}]$ and such that the initial particle density satisfies*

$$(21) \qquad X^M(0, x) \Longrightarrow X(0, x)$$

*as $M \to \infty$ for some continuous function $X(0, x)$ with support $[-\kappa, \kappa]$. Then as $M \to \infty$,*

$$(22) \qquad X^M(t, x) \Longrightarrow X(t, x),$$

*where $X(t, x)$ is the density function of the Dawson-Watanabe process with killing rate $a$ and initial value $X(0, x)$. The convergence is relative to the topology of uniform convergence on compacts in the space $C(\mathbb{R}_+ \times \mathbb{R})$ of continuous functions.*

Since the *measure-valued* processes associated with the densities $X^M(t, x)$ are known to converge to the Dawson-Watanabe process, by Watanabe's theorem, to prove Theorem 6 it suffices to establish tightness. This is done by a somewhat technical application of the Kolmogorov-Chentsov tightness criterion, based on a careful estimation of moments. See [5] for details.

It is also possible to show that convergence of the particle density processes holds in Theorem 5.